\documentclass[11pt]{article}
\usepackage{amsmath,amssymb,epsf}
\usepackage{amsmath}
\usepackage[applemac]{inputenc}
\advance \textwidth 3cm
\advance \textheight 3cm
\usepackage{amsfonts}
\usepackage{latexsym}
\newtheorem{theoreme}{Theorem}
\newtheorem{lemme}{Lemma}

\newtheorem{definition}{D\'efinition}

\newtheorem{prop}{Property}
\newtheorem{remarque}{Remark}
\newenvironment{preuve}[1]{\par\noindent\underline{Proof #1} :\quad}%
{\unskip\nobreak\hfil\penalty50\hskip2em\null\nobreak\hfil%
$\Box$\parfillskip0pt\par\medskip}

\newcommand{\Supp}{\mathrm {Supp}}
\newcommand{\pv}{\mathrm {PV.}}

\title{ Toeplitz matrices for the study of the fractional
Laplacian on a bounded interval.}

\author{ Philippe Rambour\thanks{Universit\'{e} de Paris Sud,
      B\^atiment 425; F-91405
Orsay Cedex;
tel : 01 69 15 57 28 ; fax 01 69 15 60 19
      \mbox{e-mail : philippe.rambour@u-psud.fr}}
       \and Abdellatif Seghier\thanks{Universit\'{e} de Paris Sud,
        B\^atiment 425; F-91405
Orsay Cedex;
tel : 01 69 15 60 09 ; fax 01 69 15 72 34
       \mbox{ e-mail : abdellatif.seghier@wanadoo.fr}}}
\begin{document}
\maketitle
          \renewcommand{\abstractname}{Abstract}
          \begin{abstract}
\textbf{ Toeplitz matrices for the study of the fractional Laplacian on a bounded interval, an application to fractional equations.}\\
In this work we get a deep link between  $(-\Delta)^{\alpha}_{]0,1[}$ the fractional Laplacian  on the interval $]0,1[$ and $T_{N}(\varphi_{\alpha})$ the Toeplitz matrices of symbol $\varphi_{\alpha} : \theta \mapsto \vert 1-e^{i\theta}\vert ^{2\alpha}$ when $N$ goes to the infinity and for $\alpha \in ]0, \frac{1}{2}[ \cup ]\frac{1}{2},1[$. In the second part of the paper we provide a Green function for the fractional equation 
$\left((-\Delta)^{\alpha} _{]0,1[}\right)(\psi) =f$  for $\alpha\in ]0,\frac{1}{2}[$ and $f$ a sufficiently smooth function on $[0,1]$. The interest is that this Green's function is the same as the Laplacian opertaor of order $2n, n \in \mathbb N$.
 \end{abstract}

   \vspace{1cm}


\vspace{1cm}


equations, Fractional differential operators.


\section{Introduction and statement of the main results}
For $\alpha \in ]0, \frac{1}{2}[\cup ]\frac{1}{2},1[$  we recall the definition of the one-dimensional fractional Laplace operator $(-\Delta)^{\alpha}$ (with $\Delta = 
-\frac{d^{2}}{dx^{2}}$). It is defined pointwise by the principal value of the integral,
if convergent, (see \cite{Kw1}) 
$$(-\Delta)^{\alpha} (u)(x) = C_{1}(\alpha) 
\pv \int_{\mathbb R}
\frac{u(x) - u(y)}{\vert x-y\vert^{1+2\alpha} } dy\quad \mathrm{with} 
\quad C_{1}(\alpha) =\frac{2^{2\alpha} \Gamma(\frac{1+2\alpha}{2})}
{\sqrt \pi \vert \Gamma(-\alpha)\vert }  \quad x\in \mathbb R. $$ 
$(-\Delta)^{\alpha} (u)(x) $ is convergent if, for instance, $f$ is smooth in a neigborhood of $x$ and bounded on $\mathbb R$.\\
 More generally we can refer to \cite{Kw2} 
for the different equivalent definitions of the fractional Laplace operator on the real line. Up a constant, this operator is the left inverse of the 
Riesz operator on the real line, often denoted by 
 $I^{-2\alpha} (\alpha\in ]0, \frac{1}{2}[),$ and defined by  
$ I^{-2\alpha} (\psi) (x) = \frac{1}{2 \Gamma (2\alpha) \cos (\alpha\pi)}
\int_{-\infty}^{+\infty} \frac{\psi (t)}{\vert t-x\vert ^{1+2\alpha}}dt$, for $x\in \mathbb R$ and 
$\psi \in L^{p} (\mathbb R)$, with $1\le p<\frac{1}{\alpha}$ (see \cite{Sam1}). 
 \\
For $f\in \mathcal C^{\infty} _{c} (\mathbb ]0,1[)$, the fractional Laplacian operator 
on $[0,1]$, denoted by $(-\Delta)^{\alpha}_{]0,1[}$, is defined to be the restriction of 
$ (-\Delta)^{\alpha}(f)$ to $]0,1[$ (see \cite{Kw1}). Again $(-\Delta)_{]0,1[}^{\alpha}$ 
extends to an unbounded operator on $L^{2} (]0,1[)$.
We can easily obtain that for $f\in  \mathcal C^{\infty} _{c} (\mathbb ]0,1[)$ and for $x \in ]0,1[$ 
\begin{equation}\label{LAPLACE}
 \left((-\Delta)_{]0,1[}^{\alpha}(f) \right) (x) = C_{1} (\alpha) \left( \pv \int_{0}^{1}
\frac{f(x) - f(y)}{\vert x-y\vert^{1+2\alpha} } dy + \frac{f(x)}{2\alpha} \left (x^{-2\alpha}
+(1-x)^{-2\alpha}\right)\right).
\end{equation}
More generally for $f$ a function defined on $]0,1[$ and $x\in ]0,1[$ we denote by $\left((-\Delta)_{]0,1[}^{\alpha}(f) \right) (x)$
the quantity (\ref{LAPLACE}) if it exists.
In the first part of the article we get a link between 
 this operator and  the Toeplitz matrices of order $N+1$ 
with symbol $\varphi_{\alpha} = \vert 1-\chi\vert ^{2\alpha} $ where $\chi$ is defined on  $[0,2\pi[ $ by  
$\theta \mapsto e^{i\theta}$. We recall
 that a Toeplitz matrix of order $N$ with symbol $h\in L^{1}([0,2 \pi])$ is the $(N+1)\times (N+1)$ matrix $T_{N}(h)$ defined by
$\left(T_{N}(h)\right)_{i+1,j+1}= \hat h(i-j)$ where $\hat h (s)$ denote the Fourier coefficient of order $s$ of the function $h$ (see {\cite{GS},\cite{Bo.3}). 
The Toeplitz matrices of order $N \times N$ are decisive here  because they have the property to make the link between the discrete and the continuous when $N$ goes to infinity, and  they are useful to obtain a good discretization of the problem. With these tools, taking the limit at infinity, we can
obtain operators $D_{\alpha}$ for  
$\alpha\in ]-\frac{1}{2},0[\cup]0,\frac{1}{2}[\cup ]\frac{1}{2},1[$ that  can be interpreted as 
fractional derivatives. 
For a function $f$ defined on $[0,1]$, we define these operators as  the following limit (if it exists) 
  \begin{equation} \label{EQ1}\left(D_{\alpha} f\right) (x) = \lim_{N\rightarrow + \infty} N^{2\alpha} 
\sum_{l=0}^{N} \left(T_{N} (\varphi_{\alpha})\right)_{[Nx]+1,l+1} f(\frac{l}{N})\end{equation}
for $x\in ]0,1[$ and where $[Nx]= \max\{k\in \mathbb Z \vert k\le Nx \}.$\\
 For $\alpha\in ]0,\frac{1}{2}[ \cup ]\frac{1}{2},1[$ it is the framework of the theorem \ref{theo2} to show that this limit exists and is $\left((-Æ)^{\alpha}_{]0,1[}(f)\right)(x) $ for 
$f$ belonging to certain classes of functions. \\
 \begin{remarque}\label{remrem-1}
   The operator definition $D_{\alpha}$ can be easily transported over any interval $]a,b[$ as follows. If $a<b$ are two reals and $h$ is a function defined on $]a,b[$, one defines for $x$ in  $[0,1]$ $h_{a,b}(x) = h(a+(b-a) x)$.
 Then for  $\alpha \in]-\frac{1}{2},0[ \cup ]0, \frac{1}{2}[\cup ]\frac{1}{2},1[$ and $u\in ]a,b[$, we have 
 $ D_{\alpha,a,b} (h) (u) = \left( D_{\alpha} h_{a,b}\right) \left( \frac{u-a}{b-a}\right).$\\
   \end{remarque}

The second point of this paper is to invert the fractional Laplacian $(-\Delta)^{\alpha}$ on the open interval $]0,1[$ for
$\alpha\in ]0, \frac{1}{2}[$.
 For  $f$  a contracting function on $[0,1]$ 
 such that $\Supp (f) \subset ]0,1[$   we solve the equation  in $\phi$ :
\begin{equation} \label{EQ0}
(-\Delta)^{\alpha}_{]0,1[}(\phi) = f
\end{equation} 
that is also (see \cite{Blum} and\cite{Kul})
\begin{equation}\label{EQ00}
\left\{
\begin{array}{cccc}
 (-\Delta)^{\alpha}(\phi) &  =  &  f & \mathrm{in} \quad ]0,1[ \\
  \phi&  = & 0 &  \mathrm{in} \quad ]-\infty,0]  \cup [1,+\infty[.
\end{array}
\right.
\end{equation}

   We recall that  for $\mu>0$ we denote by $C^{0,\mu} (]a,b[)$ is defined as the set of the functions $\psi$ such that for all interval $[c,d]\subset ]a,b[$ with $d-c$ sufficiently small there a real $K_{[c,d]}>0$ such that the inequality  
$$ \vert \psi (x) -\psi (x') \vert \le K _{[c,d]}\vert x-x'\vert ^{\mu}$$  is true for all  $x,x'$ in $[c,d] $.\\
In all this work we say that a function $f$ is contracting on an interval $I$ if
$$ \vert f(x) - f(y)\vert \le \vert x-y\vert \quad \forall x,y \in I.$$
Finally we will say that a function is locally contracting over an interval $]a,b[$ if it belongs to $C^{0,1}(]a,b[)$.
 
If the rest of the paper we denote by $C_{\alpha}$ the constant  $ -\frac{\Gamma (2\alpha+1) \sin (\pi \alpha)}{\pi}$. We check that $C (\alpha) = - C_1(\alpha).$
Then we can write the following statements :
\begin{theoreme}\label{theo2}
 We have  :
\begin{enumerate}
\item
 if $-\frac{1}{2}< \alpha<0$ and $h\in L^{1} ([0,1])$ then for all $x\in [0,1] $ in 
 $L^{1}\left([0,1]\right)$  
$$ \left(D_{\alpha} h\right) (u) = C_{\alpha} \int_{0}^{1} \frac{h(t)}{\vert t-x\vert^{1+2\alpha} } dt.$$
\item
 For $0<\alpha<\frac{1}{2}$ and $h$  a  function in  
 $C^{0,\mu}(]0,1[)$ with $2 \alpha < \mu <1 $ 
 $$\left(D_{\alpha} (h) \right) (x) = 
\left((-\Delta)^{\alpha}_{]0,1[}(h)\right) (x),$$
uniformly in $x\in [\delta _{1},\delta _{2}]$ for $[\delta _{1},\delta _{2}] \subset ]0,1[$.
\item
 If $\frac{1}{2}<\alpha<1$ and $h \in C^{2} \left([0,1]\right)$ then $$\left(D_{\alpha} (h) \right) (x) = 
\left((-\Delta)^{\alpha}_{]0,1[}(h)\right) (x),$$
uniformly in $x\in [\delta _{1},\delta _{2}]$ for $[\delta _{1},\delta _{2}] \subset ]0,1[$.
\end{enumerate}
\end{theoreme}
\begin{theoreme}\label{theo4}
Let  $0<\alpha<\frac{1}{2}$  and let $f$ be a  real function  contracting on $[0,1]$ such that  $\Supp (f) =[a,b] \subset ]0,1[$. Then the differential equation   
$$\left((-\Delta)_{]0,1[}^{\alpha} (g)\right) = f $$ 
has only one solution locally contracting on $]0,1[$. This solution is defined  
for $z\in ]0,1[$ by
$$ g (z) =\left( \left( D_{-\alpha}(f)\right) (z) -
 \int_{0}^{1}  K_{\alpha}(z,y) f(y) dy\right),$$
where
\begin{align*}
 K_{\alpha} (u,y) &= \frac{1}{\Gamma^{2} (\alpha)} u^\alpha  y^\alpha 
 \left(\int_{1}^{+\infty} 
 \frac{(t-u)^{\alpha-1} (t-y)^{\alpha-1}}{t^{2\alpha} }dt \right.
 \\&
 + \left. \int_{0}^{+\infty}  \frac{(t+u)^{\alpha-1} (t+y)^{\alpha-1}}
 {t^{2\alpha} }dt \right). 
 \end{align*}
    \end{theoreme}
    \begin{remarque}
The consistency of the statement of the theorem \ref{theo4} ( i.e. if $\alpha\in ]0,\frac{1}{2}[$ 
the function $\left( D_{-\alpha}(f)\right) (z) -\left( K_{\alpha} (f)\right) $ is
locally contracting on $]0,1[$) is specified in the theorem demonstration).On the other hand in \cite{RQuebl} 
 we have obtained $K_{\alpha} (x,y) = O(\vert x-y\vert ^{\alpha-1})$ for 
 $0<x\neq y <1$. Hence our solution is well defined.
\end{remarque}
The theorem \ref{theo4} also gives us  that the solution of the equation 
$(-\Delta)^{\alpha} (\psi) =f$  found for the fractional Laplacian of order 
$\alpha$ defined on an interval 
 is not up a constant the Riesz operator of order $-\alpha$ on 
the same interval, there is a perturbation, unlike the result on the real line (see \cite{Sam1}).\\
Our calculation methods can also invert the Riesz operator over a bounded interval. It is an alternative to the results given for example in \cite{Sam1}.  \\  
Using the results of  \cite{RS10} we can also write the following equivalent statement.
\begin{theoreme}\label{C1}
Let  $0<\alpha<\frac{1}{2}$  and let $h$ be a  real function  contracting on $[0,1]$such that  $\Supp (f) =[a,b] \subset ]0,1[$. Then the differential equation   
$$\left((-\Delta)_{/]0,1[}^{\alpha} (g)\right) = h $$ 
has only one solution locally contracting on $]0,1[$ (hence in $C^{0}(]0,1[) $) with $g(0) =g(1) =0$. This solution is defined  by 
\begin{enumerate}
\item
 $ g (x) =  \int_{0}^{1} G_{\alpha} (x,y) h(y) dy$
 with 
 $$ G_{\alpha}(x,y) = \frac{1}{\Gamma^{2} (\alpha)} (x)^{\alpha} (y)^{\alpha} 
 \int_{\max(x,y)} ^{1}
 \frac{(t-x)^{\alpha-1}(t-y)^{\alpha-1}}{(t)^{2\alpha}} dt, \quad \mathrm{for} 
  \quad 0<x\neq y <1$$ and 
 $$  G_{\alpha}(0,0)=0.$$
\item
$g(z)=0$ for $z\le 0$ or $z\ge1$.
\end{enumerate}
 \end{theoreme}
 \begin{remarque}\label{rem-1}
 The expression of the Green kernel $G_{\alpha} $  makes it easy to verify that the solution proposed in the theorem \ref{theo4} is extendable by zero on 
 $\mathbb R\setminus ]0,1[$. This solution is well defined. In \cite{RQuebl} 
 we have obtained $G_{\alpha} (x,y) = O(\vert x-y\vert ^{2\alpha-1})$ for 
 $0<x\neq y <1$.
  \end{remarque}
  In fact the theorem \ref{C1} is a generalization of the well known case where 
 $\alpha \in \mathbb N^{\star}$ and $[a,b]= [0,1]$ (see 
 \cite{CFL,SpSt,RS04}). In this case we have a Green function 
 $G_{\alpha}(x,y) = \frac{ x^{\alpha}y^{\alpha}}{\Gamma^{2}(\alpha)}
  \int_{\max(x,y)}^{1} \frac{(t-x)^{\alpha-1} (t-y)^{\alpha-1}}{t^{2\alpha}} dt$ for $0<\max (x,y)\le 1$ and $G(0,0)=0$ such that for all function $f\in L^1[0,1]$ the function $g$ defined on 
 $[0,1]$ by $g(x) = \int_{0}^1 G_{\alpha}(x,y) f(y) dy$ is the solution of 
 the equation (\ref{EQ0}) with the bound condition $g(0)= \cdots=g^{\alpha-1}(0)=0$
  and $g(1)=\cdots = g^{\alpha-1}(1)=0.$ It is important to 
  remark that the expression of the Green function is finally the same in the case of $\alpha \in \mathbb N$ and for the case $\alpha \in ]0, \frac{1}{2}[$. \\
 To get the solution to 
the equation (\ref{EQ0}) we use the fine knowledge of the matrices  
$\left(T_{N}(\vert 1-\chi\vert ^{2\alpha})\right)^{-1}$ for $\alpha\in ]0, \frac{1}{2}[$ that we acquired in our previous works. More precisely 
the theorems \ref{theo4} and \ref{C1} may be related to the following results obtained respectively in \cite{RQuebl} and \cite{RS10},  
where we have obtained two alternative asymptotic expansions when $N$ goes to the infinity of 
$\left(T_{N}( \varphi_{\alpha}) \right )^{-1}_{k+1,l+1}$ for $k,l$ sufficiently larges and $\alpha\in ]0, \frac{1}{2}[$.
 These expressions are specified in the two following theorems.
 \begin{theoreme} \label{theoinverse1}
For $\alpha\in ]0,\frac{1}{2}[$ we have 
\begin{equation}
 \left(T_{N} (\varphi_{\alpha})\right)^{-1} _{[Nx]+1,[Ny]+1}= 
\widehat {\varphi_{-\alpha}} \left({ [Nx]-[Ny]} \right)-N^{2\alpha-1}
  K_{\alpha}(z,y) +o(N^{2\alpha-1})
 \end{equation}
uniformly in $x,y$ for $x,y \in [\delta_{1}, \delta_{2}]$.
\end{theoreme}
and 
 \begin{theoreme} \label{theoinverse2}
 For $\alpha\in ]0,\frac{1}{2}[$ we have 
\begin{equation} \label{eqq1}
\left(T_{N} (\varphi_{\alpha})\right)^{-1}_{[Nx]+1,[Ny]+1}=
N^{2\alpha-1} G_{\alpha}(x,y)+ o(N^{\alpha-1})
\end{equation} 
uniformly in $x,y$ for $0<\delta_{1}<x\neq y < \delta_{2}<1$.
   \end{theoreme}
The references \cite{Bucur1}, \cite{Sam1}, \cite{Gorenflo1} are good introductions to fractional integrals and derivatives, fractional Laplacian, and fractional differential equations. \\
 The discretization methods used here can be extended to the study of other fractional differential operators.
Thus in an other work \cite{MarchX} we found known results concerning other fractional derivatives by 
the same discretization process using an $N+1\times N+1$ Toeplitz matrix of symbol
 $h_{\alpha}= \lim_{R\rightarrow 1^{-}} h_{\alpha,R}$ whith $1>\alpha>0$ and 
 where $h_{\alpha,R}$ is the function defined by 
 $ \theta\mapsto (1-Re^{i\theta})^{\alpha}(1+R e^{-i\theta})^{\alpha},$ for $R\in ]0,1[$ and $\theta \in [0, 2 \pi[$.
For $f$ a function defined on $[0,1]$ and $0\le x\le 1$ we then study the limit 
$$ \lim_{N\rightarrow +\infty} N ^{\alpha} 
 \left(\sum_{l=0}^{N}T_{N} \left(\varphi_{\alpha}\right)_{k+1,l+1} \left(X_{N}\right)_{l}
 \right) =\left( \tilde D_{\alpha} (f)\right) (x), \quad \mathrm{with} \quad 
k=[Nx].$$
We show in \cite{MarchX} that for  $f$a locally locally contracting on $]0,1[$ this limit is
$$
\frac{2^\alpha}{\Gamma(-\alpha)}\left(
 \int_{0}^{x} \frac{ f(t)- f(x)}{\vert x-t\vert^{-\alpha-1}}dt - f(x) \left(\frac{(x)^{-\alpha}}
 {\alpha}\right)\right)
 $$
 which is nothing more than the inferior fractional Marchaud derivative of order $\alpha$ on $[0,1]$.
We can also verify that if we choose as a symbol the function 
$\overline{h_{\alpha}}$ 
 this same limit gives us the superior fractional Marchaud derivative of order $\alpha$ on $[0,1]$.
  Still in \cite{MarchX} we find, by methods similar to those used here, the inverse of these fractional derivatives.\\
 Another approach to fractional differential equations different from the classical approach can be found in \cite{PGSG} where the authors use Hankel's operators to solve on $\mathcal S^{1}$, 
the torus of dimension 1,  the equation  $i \partial _{t} u = \pi (\vert u \vert^{2} u)$ 
where  $\pi$ is the usual orthogonal projection from
$L^{2}(\mathcal S^{1})$ on the subspace $H^{2}(\mathcal S^{1})$ 
defined by $h \in H^{2}(\mathcal S^{1}) \iff \hat h (s) =0 \, \forall s<0$.\\
Integrals and fractional derivatives are currently the focus of much mathematical works.
For example, one could consult \cite{TOBIAS2,TOBIAS1,Bucur1,DIM1,Pod1,Martadelia}.
\\
 \section{Proof of the theorem  \ref{theo2}}
 \begin{enumerate}
 \item
 It is clear for $\alpha\in ]-\frac{1}{2},0[$. 
 \item 
 Now we have to obtain the result for $\alpha\in ]0, \frac{1}{2}[$.\\
Assume that   $x$ belongs to an interval 
$[\delta _{1}, \delta _{2}]$  include in $[0,1]$. We assume that $N$ is a fixed integer and we put $k=[Nx]$. In the following we denote by  $K$ the positive constant such that for  $x,y\in [\delta _{1},\delta _{2}]$ 
\begin{equation}\label{eq2}
 \vert f(x) - f(y) \vert \le K \vert x-y\vert^{\mu}.
 \end{equation}
We will also use that 
$ \hat \varphi_{\alpha} (u) = C_{\alpha} \vert u \vert ^{-2\alpha-1} \left(1+o(1)\right)$ for $u$ an  integer with a sufficiently large absolute value. 
Let $\delta=N^{-\beta }$ with $0<\beta<1$. 
We can write 
\begin{align*}
\sum_{l=0}^{N} \left(T_{N} \varphi_{\alpha}\right)_{ k+1,l+1} f ( \frac{l}{N}) &=
\sum_{l=0}^{k-[N\delta ]-1} \left(T_{N} \varphi_{\alpha}\right)_{ k+1,l+1}  f ( \frac{l}{N})
+\sum_{l=k-[N\delta ]}^{k+[N\delta ]} \left(T_{N} \varphi_{\alpha}\right)_{ k+1,l+1}  f ( \frac{l}{N})\\
&+\sum_{l+[N\delta ]+1}^{N} \left(T_{N} \varphi_{\alpha}\right)_{ k+1,l+1}  f (\frac{l}{N}).
\end{align*}
We can now observe that  
$$
\sum_{l=k-[N\delta ]}^{k+[N\delta ]} \left(T_{N} \varphi_{\alpha}\right)_{ k+1,l+1}  f ( \frac{l}{N})
= \sum_{l=k-[N\delta ]}^{k+[N\delta ]} \left(T_{N} \varphi_{\alpha}\right)_{ k+1,l+1} 
\left( f (\frac{l}{N}) -f(\frac{k}{N})+ f(\frac{k}{N})\right).$$
Using (\ref{eq2}) we have
$$
\Bigl \vert  \sum_{l=k-[N\delta ]}^{k+[N\delta ]} \left(T_{N} \varphi_{\alpha}\right)_{ k+1,l+1} 
\left( f (\frac{l}{N}) -f(\frac{k}{N})\right) \Bigl \vert
 \le 
\vert K\vert 
\Bigl \vert \sum_{l=k-[N\delta ]}^{k+[N\delta ]} \left(T_{N} \varphi_{\alpha}\right)_{ k+1,l+1} (\frac {l}{N} - \frac{k}{N})^{\mu}\Bigr \vert
$$
that implies 
$$ \Bigl \vert  \sum_{l=k-[N\delta ]}^{k+[N\delta ]} \left(T_{N} \varphi_{\alpha}\right)_{ k+1,l+1} 
\left( f (\frac{l}{N}) -f(\frac{k}{N})\right) \Bigl \vert= O(\delta^{\mu}).$$

Lastly for $\frac{2\alpha}{\mu}<\beta<1$
we obtain
$$ N^{2\alpha} \left( \sum_{l=k-[N\delta ]}^{k+[N\delta ]} \left(T_{N} \varphi_{\alpha}\right)_{ k+1,l+1} 
\left( f (\frac{l}{N}) -f(\frac{k}{N})\right)\right) = o(1).$$
On the other, since 
$\displaystyle{\sum_{n\in \mathbb Z}  \hat \varphi_{\alpha} (n) =0}$, we have clearly
\begin{align*}
\sum_{l=k-[N\delta ]}^{k+[N\delta ]} \left(T_{N} \varphi_{\alpha}\right)_{ k+1,l+1} 
 f(\frac{k}{N})&= -\sum_{l<k-[N\delta ]}\hat \varphi_{\alpha}(k-l) f(\frac{k}{N})\\
 &
- \sum_{l>k+[N\delta ]} \hat \varphi_{\alpha}(k-l)  f(\frac{k}{N})
\end{align*}
with 
$$N^{2\alpha} \sum_{l<k-[N\delta ]}\hat \varphi_{\alpha}(k-l)  f(\frac{k}{N}) = 
\frac{C_{\alpha}}{N}\left(\sum_{l=0}^{k-[N\delta ]-1} (\frac{\vert k-l\vert}{N}) ^{-2\alpha-1} f(\frac{k}{N})
- \frac{x^{-2\alpha}}{2\alpha} f(x) \right)+o(1),$$
and
$$N^{2\alpha} \sum_{l>k+[N\delta ]} \hat \varphi_{\alpha}(k-l)  f(\frac{k}{N}) = 
\frac{C_{\alpha}}{N}\left(\sum_{l=k+[N\delta ]+1} ^{N} (\frac{\vert k-l\vert}{N} )^{-2\alpha-1} f(\frac{k}{N})
- \frac{(1-x)^{-2\alpha}}{2\alpha} f(x)\right) +o(1).$$
Since $f\in C^{0,\mu}\left(]0,1[\right)$ with $\mu>2\alpha$ we know that 
$P.V \int_{0}^{1}\frac {f(t)-f(x)}
{\vert t-x\vert ^{2\alpha+1}} dt$ is convergent and 
\begin{align*}
\lim_{N\rightarrow +\infty}\frac{C_{\alpha}}{N} & \left( \sum_{l=0}^{k-[N\delta ]-1} (\frac{\vert k-l\vert}{N}) ^{-2\alpha-1} \left(f(\frac{l}{N})-f(\frac{k}{N})\right) \right.\\
+&
\left.\sum_{l=k+[N\delta ]+1} ^{N} (\frac{\vert k-l\vert}{N} )^{-2\alpha-1} 
\left(f(\frac{l}{N})-f(\frac{k}{N})\right)\right) \\
&= C_{\alpha} P.V. \int_{0}^{1}\frac {f(t)-f(x)}
{\vert t-x\vert ^{2\alpha+1}} dt 
= C_{1} (\alpha) P.V. \int_{0}^{1}\frac {f(x)-f(t)}
{\vert x-t\vert ^{2\alpha+1}} dt 
\end{align*}
that ends the proof for $\alpha\in ]0, \frac{1}{2}[$.
\item
Lastly we have to assume that $\alpha\in ]\frac{1}{2},1[$.
Always with $\delta = N^{-\gamma}, 0<\gamma<1$ we have 
\begin{align*}
N^{2\alpha} \sum_{l=0}^{N} \hat{\varphi_{\alpha}} (k-l) f(\frac{l}{N}) &=
N^{2\alpha}\left( \sum_{l=0}^{k-[N\delta] -1} \hat{\varphi_{\alpha}} (k-l) f(\frac{l}{N}) +
\sum_{l=k-[N\delta ]}^{k+[N\delta]} \hat{\varphi_{\alpha}} (k-l) f(\frac{l}{N})\right.\\
& \left. +\sum_{l=[N\delta ]+1}^{N} \hat{\varphi_{\alpha}} (k-l) f(\frac{l}{N})\right)
\end{align*}
We can write 
$$ \sum_{l=k-[N\delta ]}^{k+[N\delta]} \hat{\varphi_{\alpha}} (k-l) f(\frac{l}{N}) =
\sum_{l=k-[N\delta ]}^{k+[N\delta]} \hat{\varphi_{\alpha}} (k-l)
 \left(f(\frac{l}{N})-f(\frac{k}{N}\right) + \left(\sum_{l=k-[N\delta ]}^{k+[N\delta]} \hat{\varphi_{\alpha}} (k-l)\right) f(\frac{k}{N}).$$
 By Taylor expansion of $f$ we obtain:
\begin{align*}
 \sum_{l=k-[N\delta ]}^{k+[N\delta]} \hat{\varphi_{\alpha}} (k-l)
 \left(f(\frac{l}{N})-f(\frac{k}{N}\right) &= \sum_{l=k-[N\delta ]}^{k+[N\delta]} \hat{\varphi_{\alpha}} (k-l) \frac{l-k}{N} f'(\frac{k}{N}) \\
 &+ \sum_{l=k-[N\delta ]}^{k+[N\delta]} \hat{\varphi_{\alpha}} (k-l) 
 \left(\frac{l-k}{N}\right)^2 f^{\prime\prime}(c_{k,l,N})
 \end{align*}
 where $c_{k,l,N}\in ]\frac{l}{N}, \frac{k}{N}[$ or $c_{k,l,N}\in ]\frac{k}{N}, \frac{l}{N}[$.
 Since $\varphi^{\prime}_{\alpha} (0) =0$  and with the additional remark that $\varphi^{\prime}_{\alpha}$ is an odd function we can write
   \begin{align*} N^{2\alpha} \left( \sum_{l=k-[N\delta ]}^{k+[N\delta]} \hat{\varphi_{\alpha}} (k-l) \frac{l-k}{N} \right) &= (-i) N^{2\alpha-1} \left(\sum_{l=k-[N\delta ]}^{k+[N\delta]} \hat{\varphi^{\prime}_{\alpha}} (k-l) \right)\\
  & =i N^{2\alpha-1} \left(\sum_{l<k-[N\delta ]}\hat{\varphi^{\prime}_{\alpha}} (k-l) + \sum_{l>k+[N\delta ]} \hat{\varphi^{\prime}_{\alpha}} (k-l) \right)=0
\end{align*}
On the other hand for a good choice of $\gamma$ :
$$ N^{2\alpha} \left( \sum_{l=k-[N\delta ]}^{k+[N\delta]} \hat{\varphi_{\alpha}} (k-l) 
 \left(\frac{l-k}{N}\right)^2 f^{\prime\prime}(c_{k,l,N})\right) \le 
K' N^{2\alpha-2} (N\delta )^{2} \Vert f^{\prime\prime} \Vert_{\infty}= O(N^{2\alpha-2\gamma}) =o(1).$$
Always since $\varphi_{\alpha}(0) =0$ we obtain as for the case $0<\alpha<\frac{1}{2}$ 
\begin{align*}
N^{2\alpha} \sum_{l=0}^{N} \hat{\varphi_{\alpha}} (k-l) f(\frac{l}{N}) &=
\frac{C_{\alpha}}{N} \left( \sum_{l=0}^{k-[N\delta] -1} \hat{\varphi_{\alpha}} (k-l) 
\left(f(\frac{l}{N})- f(\frac{k}{N})\right) \right.
\\ & +  \sum_{k+[N\delta] +1} ^{N}\hat{\varphi_{\alpha}} (k-l) \left(f(\frac{l}{N})- f(\frac{k}{N})\right)\\
&- \left. \left( \frac{x^{-2\alpha}}{2\alpha} +\frac{ (1-x)^{-2\alpha} }{2\alpha} \right) f(x) \right)+o(1).
\end{align*}
Since 
$P.V \int _{0}^{1} \frac{ f(t)-f(x)} {\vert t-x\vert ^{2\alpha+1} } dt$ is  convergent for $f\in C^{2}\left([0,1]\right) $
we obtain 
\begin{align*}
N^{2\alpha} \sum_{l=0}^{N} \hat{\varphi_{\alpha}} (k-l) f(\frac{l}{N}) &=
C_{\alpha} \left( P.V \int _{0}^{1} \frac{ f(t)-f(x)} {\vert t-x\vert ^{2\alpha+1} }dt 
-  \left( \frac{x^{-2\alpha}}{2\alpha} +\frac{ (1-x)^{-2\alpha} }{2\alpha} \right) f(x) \right)\\
& = C_{1\alpha}\left( P.V \int _{0}^{1} \frac{ f(x)-f(t)} {\vert t-x\vert ^{2\alpha+1} } dt 
+  \left( \frac{x^{-2\alpha}}{2\alpha} +\frac{ (1-x)^{-2\alpha} }{2\alpha} \right) f(x) \right)
\end{align*} 
that is the expected result.
\end{enumerate}
\section{Demonstration of  the theorems \ref{theo4} and \ref{C1}}
In the following we denote by $g_{\alpha}$ the function defined for $\theta$ in $[0,2 \pi[$ by
$\theta\mapsto (1- e^{i\theta}) ^{\alpha}$, $0<\alpha<\frac{1}{2}$
 and 
$\beta^{\alpha}_{u}$  will be the Fourier coefficient 
$\widehat {g_{\alpha}^{-1}}(u)$ for $u\in \mathbb N$ . 
It is known that for a sufficiently large integer $u$ we have 
$\beta_{u}^{\alpha} = \frac{u^{\alpha-1}}{\Gamma (\alpha)} +o(u^{\alpha-1})$.
In the demonstration we use the predictor polynomial of the functions
$\varphi_{\alpha}$,
$\alpha\in ]0, \frac{1}{2}[$ and an expression of their coefficients which has been obtained in a previous work.
In the following section the reader will find some reminders about these results.
\subsection{Predictor polynomials of $\varphi_{\alpha}$}
First we have to recall the definition of a predictor polynomial of the function $f$.
\begin{definition} 
If $h$ is an integrable positive function with have only a finite number of zeros on $[0, 2\pi[$ the predictor polynomial of degree 
$M$ of $h$ is the trigonometric polynomial defined by 
$$ P_{M} = \frac{1}{ \sqrt{\left( T_{N}(h)\right)^{-1}_{(1,1)}}}
\sum_{u=0}^M \left( T_{N}(h)\right)^{-1}_{(u+1,1)} \chi^u.$$
\end{definition}
These predictor polynomials are closely related to the orthogonal polynomials 
$ \Phi_{M}, M\in \mathbb N$ with respect to the weight $h$ 
by the relation 
$\Phi_{M} (z) = z^M \overline{ P_{M} (z)}$, for 
$\vert z \vert =1$. This relation and the classic results on the orthogonal polynomial imply that $P_{M} (e^{i\theta}) \neq 0$ for all real $\theta$.
 In the proof we will also need to use the fundamental property:
 \begin{theoreme}\label{predi}
 If $P_{M}$ is the predictor polynomial of a function $h$ then
$$ \widehat {\frac{1}{\vert P_{M}\vert ^2}} (s)
= \hat h (s) \quad \forall s \quad -M\le s\le M.$$
\end{theoreme}
That provides 
$$ T_{M}\left( \frac{1}{\vert P_{M}\vert ^2}\right) 
= T_{M}(h).$$
 \cite{Ld} is a good reference about the predictor polynomials. 
 In the proof of the lemma \ref{lemme2} we  use the coefficients of the predictor polynomials of the functions $\varphi_{\alpha}$. 
The expression we need and use is  an exact expression and has been obtained in\cite{RS10}.
Thank of the results of this last paper we can write, 
\begin{equation}\label{eqF2} \forall k, l \in [0,N] 
\quad ( T_{N} \varphi_{\alpha})^{-1}_{k+1,1} =
\beta_{k}^{\alpha}  -\frac{1}{N} 
\sum_{u=0}^{k} \beta^{\alpha}_{k-u} 
F_{\alpha,N} (\frac{u}{N}).
\end{equation}
For $z\in [0,1]$  the quantity $F_{\alpha,N}(z) $ is defined by  
\begin{equation}\label{eqF3} 
F_{\alpha,N} (z) = \sum_{m=0}^{+\infty} F_{m,N,\alpha} (z) 
\left( \frac{\sin \pi \alpha}{\pi}\right) ^{2m+2} 
\end{equation}
where
\begin{align*} 
F_{m,N,\alpha}(z) &= \sum_{w_{0}=0}^{+\infty}
\frac{1}{1+w_{0}+\frac{1+\alpha}{N}}
\sum_{w_{1}=0}^{+\infty} \frac{1}{w_{0}+w_{1}+N+1+\alpha}\times\cdots\\
&\cdots \sum_{w_{2m}-1=0}^{+\infty} \frac{1}{w_{2m-2}+w_{2m-1}+N+1+\alpha}\\
&\sum_{w_{2m}=0}^{+\infty} \frac{1}{w_{2m-1}+w_{2m}+N+1+\alpha}
\frac{1}{1+\frac{w_{2m}}{N}+\frac{1+\alpha}{N} -z}.
\end{align*}
Using integrals we can bounded  these sums by
\begin{align*}
\tilde F_{m,N,\alpha}(z,z') &=  \int_{0}^{+\infty} \frac{1}{1+t_{0}+\frac{1+\alpha}{N}}
\int_{0}^{+\infty} \frac{1}{1+t_{0}+t_{1}} \times \cdots\\
&\int_{0}^{+\infty}  \frac{1}{1+t_{2m-2}+t_{2m-1}}
 \int_{0}^{+\infty}\frac{1}{1+t_{2m-1}+t_{2m}}  \frac{1}{1+t_{2m}+z} dt_{0}dt_{1}
\cdots dt_{2m-1} dt_{2m}.
\end{align*}
and we obtain the following upper bound, that we use 
the demonstration 
\begin{equation} \label{eqF}  
\forall z \in[0,1] \quad \vert F_{N,\alpha} (z)\vert \le
 K_{0} \left ( 1+ \Bigl\vert \ln\left( 1-z+\frac{1+\alpha}{N}\right)\Bigr\vert \right) .
 \end{equation}
 

\subsection{Existence of the solution}
Let us recall the following formula which is an adaptation of the Gohberg-Semencul formula.
\begin{prop} \label{prop2}
Let $K_{N}= \displaystyle{\sum_{u=0}^N \omega_{u}\chi^u }$ be a trigonometric polynomial of degree $N$  such that 
$K_{N}(e^{i \theta})\neq 0$ for all $\theta \in \mathbb R$.
We have, for $0\le k \le l \le N$ 
$$ T_{N}\left ( \frac{1}{\vert K_{N}\vert ^2}\right)^{-1}_{k+1,l+1}
=
\sum_{u=0}^k \omega_{k-u}\bar \omega_{l-u}
- \sum_{u=0}^k \omega_{u+N-l}\bar \omega_{u+N-k}.$$
\end{prop}
  In the next of the proof we apply this property to
  $P_{N,\alpha}$ the predictor polynomial of $\varphi_{\alpha}$,
  and we denote by $ \gamma_{u}^{\alpha}$ the 
  coefficients of $P_{N,\alpha}$. \\
 The demonstration of the lemma is based on the remark that   for $k$ such that $k =[Nx], x \in ]0,1[$,
\begin{equation} \label{BASE}
\sum_{l=0}^N 
\left(T_{N}(\varphi_{\alpha})\right)_{k+1,l+1} 
\sum_{m=0}^{N} \left(T_{N}^{-1}(\varphi_{\alpha})\right)_{l+1,m+1}  f(\frac{m}{N}) = f(\frac{k}{N}). 
\end{equation}
We consider now the term 
$$ \tilde \Phi_{N}(l) = \sum_{m=0}^N \left(T_{N}^{-1}(\varphi_{\alpha})\right)_{l+1,m+1}  f(\frac{m}{N}) 
= \sum_{m=[Na]}^{[Nb]} \left(T_{N}^{-1}(\varphi_{\alpha})\right)_{l+1,m+1}  f(\frac{m}{N}).$$
Then (\ref{BASE}) can be written 
\begin{equation} \label{BASE2}
\sum_{l=0}^N 
\left(T_{N}(\varphi_{\alpha})\right)_{k+1,l+1}  \tilde \Phi_{N}(l) 
= f(\frac{k}{N}).
\end{equation}

  We have now to use the following lemma 
 \begin{lemme} \label {lemme2}
Let $\delta_{1}$, $\delta_{2}$ be two reals in $]0,1[$, and $f$ 
a locally contracting function on $]0,1[$.  Then for all 
reals $\delta_{1},\delta_{2}$ such that 
$0<\delta_{1}<\delta_{2}<1$ there is a positive constant $M$ such that 
$\frac{l}{N}, \frac{l'}{N} \in [\delta _{1},\delta _{2}]$ we have 
$$\vert \tilde \Phi_{N} (l)-\tilde \Phi_{N} (l')\vert \le M \vert \frac{l}{N} - \frac{l'}{N} \vert
  N^{2\alpha}.$$
uniformly in $l, l'$ in $[N\delta _{1},N\delta _{2}]$
\end{lemme}
This lemma is shown in the appendix of this article.
By the same methods as in the proof of the theorem  \ref{theo2} for $\alpha \in ]0, \frac{1}{2}[$ we get, for a
good choice of the real $\delta= N^{-\beta}$ 
\begin{align*} 
\sum_{l=0}^{N} \hat \varphi_{\alpha} (k-l) 
\tilde \Phi_{N} (l) &= 
-C_{\alpha} \sum_{l=0}^{k-N\delta -1} \vert k-l\vert ^{-2\alpha-1} \left( \tilde\Phi_{N} (l)-\tilde\Phi_{N} (k)\right))\\
&-C_{\alpha} \sum_{k+N\delta +1}^{N} \vert k-l\vert ^{-2\alpha-1} \left( \tilde\Phi_{N} (l)-\tilde\Phi_{N} (k)\right)\\
&+ C_{\alpha} \tilde\Phi_{N} (k) \left( \sum_{l<0} \vert k-l\vert ^{-2\alpha-1} + \sum_{l>N} \vert k-l\vert ^{-2\alpha-1} \right)
+o(1).
 \end{align*}
 For $y\in [0,1]$
 we put 
$$ \tilde \Phi (y) =
 \int_{0}^1 H_{\alpha}(y,z) f(z)dz$$
 with $H_{\alpha}(u,v) = C_{-\alpha} \vert u-v\vert ^{2\alpha-1} -K_{\alpha}(u,v)$ 
for $ u,v\in ]0,1[$. 
 We have now to prove the following property 
  \begin{prop}\label{P2}
  \begin{enumerate}
  \item
    For a good choice of $\delta$ and for $k=[Nx]$, $0<x<1$ we have, 
  \begin{enumerate}
  \item
  $$ \sum_{l=0}^{k-[N\delta]} \vert k-l\vert ^{-2\alpha-1}
 \left ( \tilde \Phi_{N}(l) - \tilde \Phi_{N} (k)\right) = 
 \int_{0}^{x-\delta} \frac{ \tilde \Phi (y) -\tilde \Phi (x)}
 {\vert x-y\vert ^{2\alpha+1}} dy +o(1).$$
 \item
 $$ \sum_{k+[N\delta]}^N \vert k-l\vert ^{-2\alpha-1}
 \left ( \tilde \Phi_{N}(l) - \tilde \Phi_{N} (k)\right) = 
 \int_{x+\delta}^{1} \frac{ \tilde \Phi (y) -\tilde \Phi (x)}
 {\vert x-y\vert ^{2\alpha+1}} dy +o(1).$$
 \end{enumerate}
 \item
 For $k=[Nx]$ $0<x<1$ we can write uniformly in $x\in [\delta_{1}, \delta_{2}]\subset ]0,1[$, 
 \begin{align*}
 \tilde \Phi_{N}(k) & \left( \sum_{l<0} 
 \vert k-l\vert ^{-2\alpha-1} + \sum_{l<0} 
 \vert k-l\vert ^{-2\alpha-1}\right) \\
 &= \tilde \Phi (x) \left( \int_{-\infty }^0 \vert x-y\vert ^{-2\alpha-1}
 dy + \int_{1}^{+\infty} \vert x-y\vert ^{-2\alpha-1} dy\right)+o(1)
 \end{align*} 
 \end{enumerate}
 \end{prop}

We have now to prove the property \ref{P2}.
   \begin{preuve}{of the property \ref{P2}}
   First we consider an interval $[\delta_{1}, 1-\delta_{1}]$ which 
   contains $[a,b] = \Supp f$ and such that $k \in [\delta_{1}, 1-\delta_{1}]$. For bound $ \Bigr \vert \left(T_{N} ^{-1} (\varphi_{\alpha})\right)_{l+1,m+1}\Bigl \vert $ when 
   $l \le N \delta_{1}$ and  
   $m \in [Na,Nb]$ 
   we will use the following lemma that we have demonstrated in \cite{RQuebl}.
 \begin{lemme}\label{Q1}
 Let $\delta_{0}$ be a positive real and $\alpha \in ]0, \frac{1}{2}[$.Then we have a constant   
 $K_{1,\alpha}$ depending only from $\alpha$ such that, for a sufficiently large  $N$  $$ \left(T_{N}^{-1} (\varphi_{\alpha})\right)_{k+1,l+1} \le K_{1,\alpha}
 \vert l-k\vert^{\alpha-1} N^\alpha \delta^{\alpha/2}$$
 for all  $(k,l)\in \mathbb N^{2}$ with
 $0\le \min (k,l)<N\delta$ and 
 $2N\delta_{0}<\max (k,l)<N-2N\delta_{0}$.
 \end{lemme}
This lemma gives us, for $0\le l \le N \delta_{1}$ and 
$ m \in [Na,Nb]$ 
$$ \Bigr \vert \left(T_{N}  (\varphi_{\alpha})\right)^{-1}_{l+1,m+1}\Bigl \vert 
   = O\left( N^{2\alpha-1} \delta_{1}^{\alpha/2})\right),$$
   and, by the symmetries of the matrix 
   $\left(T_{N}  (\varphi_{\alpha})\right)^{-1}$ we have also, for 
   $N(1-\delta_{1})< l \le N$ and $m\in[Na, Nb]$,
   $$ \Bigr \vert \left(T_{N}  (\varphi_{\alpha})\right)^{-1}_{l+1,m+1}\Bigl \vert 
   = O\left( N^{2\alpha-1} \delta_{1}^{\alpha/2})\right),$$
   On the other hand 
   the theorems \ref{theoinverse1}  or \ref{theoinverse2} provide, for all $m \in [Na,Nb]$,  
   $$ \Bigr \vert \left(T_{N} ^{-1} (\varphi_{\alpha})\right)_{k+1,m+1}\Bigl \vert 
   = O\left( N^{2\alpha-1}\right).$$
   Hence we can say that 
  for $0\le l\le N\delta_{1}\le N $ 
   or $ N(1-\delta_{1})\le l \le  N $ we have 
   $\vert \tilde \Phi_{N}(l)\vert \le 
   O(N^{2\alpha }\delta _{1}^{\alpha/2}). $\\
   Then 
   \begin{enumerate}
   \item
   $$\Bigl \vert  \sum_{l=0} ^{N\delta _{1}} \vert k-l\vert ^{-2\alpha-1}
   \tilde \Phi_{N} (l)\Bigr \vert\le O(\delta _{1}^{3\alpha/2})
   \quad \mathrm {and} \quad 
   \Bigl \vert  \sum_{l=N(1-\delta_{1})} ^{N} \vert k-l\vert ^{-2\alpha-1}
   \tilde \Phi_{N} (l)\Bigr \vert\le O(\delta _{1}^{3\alpha/2}).$$
   \item
$$\Bigl \vert  \sum_{l=0} ^{N\delta _{1}} \vert k-l\vert ^{-2\alpha-1}
   \tilde \Phi_{N} (k)\Bigr \vert\le O(\delta _{1}).$$
   \end{enumerate}
Hence we can write, for $\delta_{1}$ sufficiently small,
$$ \sum_{l=0} ^{k-N\delta} \vert k-l\vert ^{-2\alpha-1}
   \left( \tilde \Phi_{N} (l) - \tilde \Phi_{N} (k)\right) = 
\sum_{l=N\delta _{1}+1} ^{k-N\delta} \vert k-l\vert ^{-2\alpha-1}
   \left( \tilde \Phi_{N} (l) - \tilde \Phi_{N} (k)\right) 
    +O\left(\delta _{1} ^{\alpha/2}\right)$$
    and  
    $$ \sum_{l=k+N\delta}^N \vert k-l\vert ^{-2\alpha-1}
   \left( \tilde \Phi_{N} (l) - \tilde \Phi_{N} (k)\right) = 
\sum_{l=k+N\delta}^{N (1-\delta_{1})} 
\vert k-l\vert ^{-2\alpha-1}
   \left( \tilde \Phi_{N} (l) - \tilde \Phi_{N} (k)\right) 
    +O\left(\delta _{1}^{\alpha/2} \right)$$ 
    With the uniformity in the theorem \ref{theoinverse1} 
    we obtain 
   \begin{align*}
   & \sum_{l=N\delta _{1}+1} ^{k-N\delta} \vert k-l\vert ^{-2\alpha-1}
   \left( \tilde \Phi_{N} (l) - \tilde \Phi_{N} (k)\right) \\
   &= 
   \sum_{l=N\delta _{1}+1} ^{k-N\delta} \vert k-l\vert ^{-2\alpha-1}
   N^{2\alpha}\left( \int_{a}^{b} \left(H_{\alpha} (t,\frac{l}{N})- H_{\alpha}(t,\frac{k}{N}) \right) f(t) dt\right) +o(1)
   \end{align*}
   and 
   \begin{align*}
   & \sum_{l=k+N\delta} ^{N(1-\delta_{1})}\vert k-l\vert ^{-2\alpha-1}
   \left( \tilde \Phi_{N} (l) - \tilde \Phi_{N} (k)\right) \\
   &= 
   \sum_{l=k+N\delta} ^{N(1-\delta_{1})}
   \vert k-l\vert ^{-2\alpha-1}
   N^{2\alpha} \left(\int_{a}^{b} \left(H_{\alpha} (t,\frac{l}{N})- H_{\alpha}(t,\frac{k}{N}) \right) f(t) dt\right) +o(1)
   \end{align*}
   To conclude we obtain 
\begin{align*}  & \sum_{l=N\delta _{1}+1} ^{k-N\delta} \vert k-l\vert ^{-2\alpha-1}
   N^{2\alpha}\int_{a}^{b} \left(H_{\alpha} (t,\frac{l}{N})- H_{\alpha}(t,\frac{k}{N}) \right) f(t) dt =\\
 &  \int_{\delta _{1}} ^{x-\delta} (x-y)^{-2\alpha-1} \int_{a}^{b} \left(H_{\alpha}(t,y) -H_{\alpha}(t,x) \right) f(t) dt dy+o(1).\end{align*}
that is also,
 \begin{align*}   &\sum_{l=N\delta _{1}+1} ^{k-N\delta} \vert k-l\vert ^{-2\alpha-1}
   N^{2\alpha}\int_{a}^{b} \left(H_{\alpha} (t,\frac{l}{N})- H_{\alpha}(t,\frac{k}{N}) \right) f(t) dt =\\
 &  \int_{\delta _{1}} ^{x-\delta} (x-y)^{-2\alpha-1} \int_{0}^{1} \left(H_{\alpha}(t,y) -H_{\alpha}(t,x)\right) f(t) dt dy+o(1).\end{align*}
  Since the function $t \to H(t,z)$ is in $L^1 [0,1]$ for all $z$ 
  in $[0,1]$  the beginning of the proof implies that for 
  $\delta_{1}\to 0 $, 
\begin{align*}  & \sum_{l=0} ^{k-N\delta} \vert k-l\vert ^{-2\alpha-1}
   \left( \tilde \Phi_{N} (l) - \tilde \Phi_{N} (k)\right) =\\
   & 
   \int_{0} ^{x-\delta} (x-y)^{-2\alpha-1} \int_{0}^{1} \left(H_{\alpha}(t,y) -H_{\alpha}(t,x)\right) f(t) dt dy+o(1).\end{align*}
   and identically 
   \begin{align*}  & \sum_{l=k+N\delta}^N \vert k-l\vert ^{-2\alpha-1}
   \left( \tilde \Phi_{N} (l) - \tilde \Phi_{N} (k)\right) =\\
   &
   \int_{x+\delta} ^1 (x-y)^{-2\alpha-1} \int_{0}^{1} \left( H_{\alpha}(t,y) -H_{\alpha}(t,x)\right) f(t) dt  dy+o(1).\end{align*}
   The second point of the property is clear.
   \end{preuve}
   It is now a direct consequence of  the property (\ref{P2}) and of the equality \ref{TOUT} is that for $x\in ]0,1[$
  \begin{equation} \label{FINALE-1} 
  f(x) = -C_{\alpha}
  \left( P.V. \int_{0}^1
   \frac{\tilde \Phi (y) -\tilde \Phi(x)}{\vert x-y\vert ^{2\alpha-1}}dy
   - \frac{\tilde \Phi (x) }{-2\alpha}
    \left ( x^{-2\alpha}+(1-x)^{-2\alpha}\right)\right). 
    \end{equation}
   A direct consequence of the lemma \ref{lemme2} and is that the function 
    $\tilde \Phi$ is locally contracting the equation (\ref{FINALE-1}) implies 
    \begin{equation}\label{FINALE} 
    f(x) =\left( (-Ê\Delta)^\alpha (\tilde \Phi )\right) (x). 
    \end{equation}
To ends the demonstration we have to notice that for 
$y \in [0,1]$ we have the following relation 
$ \tilde \Phi (y) = C_{-\alpha} \int_{0}^1 
\frac{f(z)} {\vert y-z\vert ^{-2\alpha+1} }dz - \int_{0}^1
K_{\alpha} (y,z) dz$
that is also 
$ \tilde \Phi (y)= \left(D_{-\alpha} (f)\right)(y) 
- \int_{0}^1
K_{\alpha} (y,z) dz$.

 \subsection{Proof of the theorem \ref{C1}}
  Now we have to involve the generalized Green kernel $G_{\alpha}$ to obtain the theorem (\ref{C1}).
For $x,y \in [\delta_{1}, \delta_{2}] \subset ]0,1[$ the theorems \ref{theoinverse1}   
provides the equality
$$ N^{2\alpha-1}  G_{\alpha}(x,y) f(y)+o(N^{2\alpha-1}) = N^{2\alpha-1} C_{-\alpha}\vert x-y\vert^{2\alpha-1}f(y) 
- N^{2\alpha-1}  K_{\alpha}(x,y) f(y)$$
for  $x,y \in [\delta_{1}, \delta_{2}] \subset ]0,1[$, $x\neq y$. 
That means
$$ G_{\alpha}(x,y) f(y)+o(1) = \vert x-y\vert^{2\alpha-1}f(y) 
-  K_{\alpha}(x,y) f(y).$$
Since for a fixed $x$ the functions $y\mapsto G(x,y)$, and $y \mapsto K_{\alpha} (x,y)$ are in $L^{1}([0,1])$  we can write, for  $x,y \in [\delta_{1}, \delta_{2}] \subset ]0,1[$, $x\neq y$, 
$$\int_{0}^{1} G_{\alpha} (x,y) f(y) dy = \int_{0}^{1}C_{-\alpha}\frac{f(y)} {\vert x-y\vert ^{1-2\alpha}} dy 
- \int_{0}^{1} K_{\alpha}(x,y) f(y) dy$$
That is also
 $$\int_{0}^{1} G_{\alpha} (x,y) f(y) dy = \left( D_{\alpha} (f)\right) (x) - \int_{0}^{1} K_{\alpha}(x,y) f(y) dy.$$
That ends the proof.
\subsection{Unicity of the solution}
Let $f$ and $\phi$ be two functions defined on  $[0,1]$, with $\phi$ a locally contracting on $]0,1[$ such that  $\Phi (0)= \Phi (1)=0$.
Moreover we assume than $f$ is a  contracting function on $[0,1]$ with $\Supp (f) = [a,b]\subset ]0,1[$.
We make the hypothese $ (-\Delta)^{\alpha} (\phi) = f$. The aim of this part is to prove 
that the function $\phi$ checks the theorem assumptions on $]0,1[$. \\
For fixed $N$ we define the vector  $X_{N}$ (reps. $Y_{N}$) of length $N+1$ by   
$\left(X_{N}\right)_{k}= \phi (\frac{k}{N})$  (resp. $\left(Y_{N}\right)_{k}=
 f(\frac{k}{N})$, $0\le k \le N$.\\
For a sufficiently large $N$ we can write 
$N^{2\alpha } \left( T_{N}(\varphi_{\alpha}) (X_{N})\right)_{k} = (Y_{N})_{k}+ (R_{N})_{k}$,
that gives us, for $0\le m\le N$
$$ (X_{N})_{m}= N^{-2\alpha} \left( T_{N}(\varphi_{\alpha}) ^{-1}(Y_{N})\right)_{m} +
\left( T_{N}(\varphi_{\alpha})^{-1} (R_{N})\right)_{m}.$$
First we have to compute $N^{2\alpha } \left( T_{N}(\varphi_{\alpha}) (X_{N})\right)_{k} $ to evaluate precisely the order of  $\vert R_{N}(k)\vert $.\\
 To do this we have to recall that for all positive real $\epsilon$ 
 we have an integer $M_{0}$ such that for $\vert u \vert \ge M_{0}$ $\hat \varphi_{\alpha} (u) = C_{\alpha} u^{-2\alpha-1}+C'_{\alpha} u^{-2\alpha-2}
 (1+r_{u})$ with $\vert r_{u}\vert \le \epsilon$ (\cite {Zyg2}).
 In the following of the proof we denote by $M$
 an integer  $M= N \delta$ with 
 $\delta = N ^{-\beta}$, $2 \alpha<\beta<1$, these conditions being set to be consistent with 
 the demonstration of the theorem \ref{theo2}. 
 Now let  $k$ such that  
 $\frac{k}{N}\in [\delta_{1}, 1-\delta_{1}]$, with
 $0<\delta_{1}<1-\delta_{1}<1$.
 According to the proof of the theorem \ref{theo2} we have to consider the seven following quantities.
 
\begin{enumerate}
\item
The sum 
$$N^{2\alpha} \left( \sum_{l=k-M} ^{k+M} 
(T_{N}\varphi_{\alpha})_{k+1,l+1}\left( \phi(\frac{l}{N}) -\phi(\frac{k}{N})\right)
\right)$$ which is $O(N^{-\beta+2\alpha})=o(1)$.
\item
The difference  
$$ N^{2\alpha} \sum_{l=0}^{k-M-1} (T_{N}\varphi_{\alpha})_{k+1,l+1} 
\left(\phi(\frac{l}{N}) -\phi( \frac{k}{N} ) \right)- N^{2\alpha} 
C_{\alpha}\sum_{l=0} ^{k-M-1} \vert k-l\vert ^{-2\alpha-1} 
\left(\phi(\frac{l}{N}) -\phi( \frac{k}{N}) \right)
$$
which is bounded by  \\
$N^{2\alpha}\displaystyle{ \sum_{l=0}^{k-M-1} O(\vert k-l\vert ^{-2\alpha-2})} \left(\phi(\frac{l}{N}) -\phi( \frac{k}{N}) \right)
=O\left (N^{-1} \int_{0}^x \frac{\Phi (t)-\Phi(x)}
{\vert t-x\vert ^{2\alpha+1}}\right)= o(1)$.
\item
Identically 
\begin{align*}
N^{2\alpha} \sum_{l=k+M}^{N} (T_{N}\varphi_{\alpha})_{k+1,l+1} 
\left(\phi(\frac{l}{N}) -\phi(\frac{k}{N}) \right)
&- N^{2\alpha} C_{\alpha}\sum_{l=k+M} ^{N} \vert k-l\vert ^{-2\alpha-1}
\left(\phi(\frac{l}{N}) -\phi( \frac{k}{N})  \right)\\
&= O(N^{-1})= o(1).
\end{align*}
\item
The difference between the quantity  
$$ N^{2\alpha}C_{\alpha}\sum_{0\le l\le k-M, N\ge l\ge k+M} \vert k-l\vert ^{-2\alpha-1} 
\left(\phi(\frac{l}{N}) -\phi( \frac{k}{N})  \right)$$
and the quantity
$$N^{2\alpha} \left( \int_{0}^{k-M} \vert k-l\vert ^{-2\alpha-1} 
\left(\phi(\frac{l}{N}) - \phi( \frac{k}{N}) \right) dl + 
\int_{k+M}^N 
\vert k-l\vert ^{-2\alpha-1} 
\left(\phi(\frac{l}{N}) -\phi( \frac{k}{N})  \right)dl\right) $$
is in $
O(  k^{-2\alpha-1} \vert \Phi(x)\vert +N^{-1+2\alpha \beta}) =O(N^{-2\alpha-1})= o(1)$ (residual term of the Euler 
and Mac-Laurin formula).
\item
The difference between 
$\displaystyle{N^{2\alpha} \sum_{l<0} \hat \varphi _{\alpha} 
(k-l) f(\frac{k}{N})}$
and $\displaystyle{N^{2\alpha} \sum_{l<0} \vert k-l\vert ^{-2\alpha-1} f(\frac{k}{N})}$ is in $O(N^{-1})$. Idem for the difference between 
$\displaystyle{N^{2\alpha} \sum_{l>N} \hat \varphi _{\alpha} (k-l) f(\frac{k}{N})}$
and $\displaystyle{N^{2\alpha} \sum_{l>N} \vert k-l\vert ^{-2\alpha-1} f(\frac{k}{N})}$
\item
The difference between 
$\displaystyle{N^{2\alpha} \sum_{l<0} \vert k-l\vert ^{-2\alpha-1} f(\frac{k}{N})}$ and $\int_{\infty}^0 \vert x-t\vert ^{-2\alpha-1} f(x) dt$
is in $O(N^{-1})$ (residual term of the Euler Mac-Laurin formula). Idem for $\displaystyle{N^{2\alpha} \sum_{l>N} \vert k-l\vert ^{-2\alpha-1} f(\frac{k}{N})}$ and $\int_{1}^{+\infty} \vert x-t\vert ^{-2\alpha-1} f(x) dt$
\item
Finally the difference between the integrals of the previous point and 
$ \int_{0}^1 \frac{\phi(t)-\phi(x)}{\vert x-t\vert ^{2\alpha+1}}
dt $ is  $O \left(N^{\beta(2\alpha-1)}\right)=o(1) $.
\end{enumerate}
Thanks to these results we can affirm the property
\begin{prop} \label{P0}
for all integer $k$, 
$M\le k \le N-M$ $R_{N}(k)= o(1) $, uniformly in $k$.
\end{prop}
On the other hand if  $0\le k \le N\delta_{1}$ or $N-N\delta_{1}\le k\le N$ 
$$R_{N}(k) = f(\frac{k}{N}) - \sum_{l=0}^N \left(T_{N}(\varphi_{\alpha}) \right)_{k+1,l+1}
\phi (\frac{l}{N})$$
which is bounded by a constant $K_{0}$ no depending from $N$.\\
Now we have to study  
$$S_{m} = N^{-2\alpha}\sum_{k=0}^{N} \left( T_{N}^{-1} (\varphi_{\alpha})\right)^{-1}_{m+1,k+1}
 (R_{N})_{k} $$
for $m\in \mathbb N$ such that $\displaystyle{ \lim_{N\rightarrow }\frac{m}{N}} = x$ and $x\in ]0,1[$. 
 To evaluate $S_{m}$  we will use the lemma \ref{Q1} the property \ref{P0} and the theorem \ref{theoinverse1}.
 Then we split the sum $S_{m}$ as follows  
\begin{align*} 
S_{m} &=N^{-2\alpha} \left( \sum_{k=0}^{N\delta_{2}-1} \left(T_{N}(\varphi_{\alpha})\right)^{-1}_{m+1,k+1} R_{N}(k)\right.\\
&+ \left.\sum_{k=N\delta_{2}}^{N-N\delta_{2}} \left(T_{N}(\varphi_{\alpha})\right)^{-1}_{m+1,k+1} R_{N}(k)+\sum_{k=N-N\delta_{2}+1}^{N} \left(T_{N}(\varphi_{\alpha})\right)^{-1}_{m+1,k+1} R_{N}(k)\right)
\end{align*}
where $0<\delta_{2}<\frac{m}{N}< 1-\delta_{2}<1$.\\
With the theorem \ref{theoinverse1} and  the property 
\ref{P0}  we get 
$$ N^{-2\alpha}\displaystyle{\sum_{k=N\delta_{2}}^{N-N\delta_{2}} \left(T_{N}(\varphi_{\alpha})\right)^{-1}_{m+1,k+1} R_{N}(k)}=O(N^{-\tau} ) =o(1).$$\\
Then the lemma \ref{Q1} and the symmetries of the matrix $ \left(T_{N}(\varphi_{\alpha})\right)^{-1}$ provides us that the sums
$N^{-2\alpha}\displaystyle{\sum_{k=0}^{N\delta_{2}-1} \left(T_{N}(\varphi_{\alpha})\right)^{-1}_{m+1,k+1} R_{N}(k)},$ 
and \\
$N^{-2\alpha}\displaystyle{\sum_{k=N-\delta_{2}+1}^{N} \left(T_{N}(\varphi_{\alpha})\right)^{-1}_{m+1,k+1} R_{N}(k)},$ are respectively of same order that the quantities \\
$N^{-2\alpha} \displaystyle{\sum_{k=0}^{N\delta_{1}-1} \vert m-k\vert ^{\alpha-1}
N^{\alpha} \delta_{2}^{\alpha/2}}=O(\delta_{1}^{\alpha/2})$  and
$N^{-2\alpha}\displaystyle{\sum_{k=N-N\delta_{2}+1}^{N}  \vert m-k\vert ^{\alpha-1}
N^{\alpha} \delta_{1}^{\alpha/2}}=O(\delta_{1}^{\alpha/2})$.\\
Hence for $\delta_{2}$  sufficiently close to zero
we can  bound $S_{m}$ by  $\epsilon\to 0 $. That means that, for considered $m$ such that 
$\lim_{N\to \infty} \frac{m}{N} =x, x\in ]0,1[$ we have 
$$\phi(x) = \lim_{N\to +\infty} N^{-2\alpha}
\left( \left(T_{N}(\varphi_{\alpha})\right)^{-1} Y_{N}\right)_{m}=\lim_{N\to \infty} 
N^{-2\alpha} \sum_{k=0} ^{N} \left(T_{N}(\varphi_{\alpha})\right)^{-1}_{m+1,k+1} 
f (\frac{k}{N}).$$
And, by repeating a calculation already done, with the notations used in the proof of the existence of a solution 
$$\lim_{N\to +\infty} N^{-2\alpha} \sum_{k=0} ^{N} \left(T_{N}(\varphi_{\alpha})\right)^{-1}_{m+1,k+1} 
f (\frac{k}{N}) = \frac{1}{\Gamma^2(\alpha)} 
 \int_{0}^1 H_{\alpha}(x,z) f (z) dz. $$
That gives us the announced unicity.
\section{Appendix : proof of the lemmas \ref{lemme2}}
\subsection{Proof of the lemma \ref{lemme2}}
According to (\ref{prop2}) we have to prove that the both functions 
$l \mapsto \displaystyle{\sum_{m=0}^{N}\sum_{u=0}^{\min(l,m)} \gamma_{\min (l,m)-u}^{\alpha} \overline {\gamma_{\max(l,m)-u}}}$ and 
$l \mapsto \displaystyle{\sum_{m=0}^{N}\sum_{u=0}^{\min(l,m)} \gamma_{N-\min (l,m)-u}^{\alpha} \overline \gamma_{N-\max(l,m)-u}}$
satisfy the property of the lemma.\\
First we prove the property for the first of these functions.
To do this it is clear with  (\ref{eqF2}) that we have  to state the three
following inequalities 
\begin{equation}\label {LIP1}
\vert \psi_{1,N} (l) - \psi_{1,N} (l')\vert \le M_{1} \vert \frac{l}{N} - \frac{l'}{N} \vert N^{2\alpha}.
\end{equation}
\begin{equation}\label {LIP2}
\vert \psi_{2,N} (l) - \psi_{2,N} (l')\vert \le M_{2} \vert \frac{l}{N} - \frac{l'}{N} \vert N^{2\alpha}.
\end{equation}
\begin{equation}\label {LIP3}
\vert \psi_{3,N} (l) - \psi_{3,N} (l')\vert \le M_{3} \vert \frac{l}{N} - \frac{l'}{N} \vert N^{2\alpha}.
\end{equation}
\begin{equation}\label {LIP4}
\vert \psi_{4,N} (l) - \psi_{4,N} (l')\vert \le M_{4} \vert \frac{l}{N} - \frac{l'}{N} \vert N^{2\alpha}.
\end{equation}

with 
\begin{enumerate} 
\item
$$ \psi_{1,N} (l) = \sum_{m=0}^{N}\left( \sum_{u=0}^{\min(l,m)}
\beta^{\alpha}_{\min(l,m)-u} \overline {\beta^{\alpha}_{\max(l,m)-u}} \right)f(\frac{m}{N})$$
\item
$$ \psi_{2,N} (l) = \sum_{m=0}^{N}\left( \sum_{u=0}^{\min(l,m)}\beta^{\alpha}_{\min(l,m)-u}
\left(\overline{\frac{1}{N} \sum_{v=0}^{\max (l,m)-u}F_{N,\alpha}(\frac{v}{N} )
\beta^{\alpha}_{\max (l,m)-u-v}}\right)\right) f(\frac{m}{N}).$$
\item
$$ \psi_{3,N} (l) = \sum_{m=0}^{N}\left( \sum_{u=0}^{\min(l,m)}\left(\frac{1}{N} \sum_{v=0}^{\min (l,m)-u}F_{N,\alpha}(\frac{v}{N} )
\beta^{\alpha}_{\min(l,m)-u-v}\right)\overline{\beta^{\alpha}_{\max(l,m)-u}}
\right) f(\frac{m}{N}).$$
\item 
\begin{align*} \psi_{4,N} (l) &= \sum_{m=0}^{N}\sum_{u=0}^{\min(l,m)}\left(
\left(\frac{1}{N} \sum_{v=0}^{\min (l,m)-u}F_{N,\alpha}(\frac{v}{N} )\beta^{\alpha}_{\min(l,m)-u-v}\right)\right.\\
& \left.
\times \left(\overline{\frac{1}{N} \sum_{w=0}^{\max (l,m)-u}F_{N,\alpha}(\frac{w}{N} )\beta^{\alpha}_{\max(l,m)-u-w}}\right) \right)f(\frac{m}{N})
\end{align*}
\end{enumerate}
First we have to prove the  lemma for $\psi_{1,N}$. We can write 
$ \psi_{1,N} (l) = \psi_{1,1,N} (l)+\psi_{1,2, N} (l)$ 
with  

$$
 \psi_{1,1,N} (l) = \sum_{m=0}^{l}\left( \sum_{u=0}^{m}
\beta^{\alpha}_{m-u} \overline {\beta^{\alpha}_{l-u}} \right)f(\frac{m}{N}).$$
$$\psi_{1,2,N} (l) = \sum_{m=l}^{N}\left( \sum_{u=0}^{l}
\beta^{\alpha}_{l-u} \overline {\beta^{\alpha}_{m-u}} \right)f(\frac{m}{N}).$$
Let be $l'\le l $ two integers with $l'<l$. We have to bound the difference
$$ \vert \psi_{1,1,N} (l) -\psi_{1,1,N} (l')\vert
\quad \mathrm{for} \quad \frac{l}{N}, \frac{l'}{N} \in [\delta _{1}, \delta _{2}].$$
We can write 
\begin{align*}
 \psi_{1,1,N} (l)  &=  \sum_{m=0}^{l'}\left( \sum_{u=0}^{l'-m}
\beta^{\alpha}_{u} \overline {\beta^{\alpha}_{m+u}} \right)f(\frac{l-m}{N})
+ \sum_{m=0}^{l'}\left( \sum_{u=l'-m}^{l-m}
\beta^{\alpha}_{u} \overline {\beta^{\alpha}_{m+u}} \right)f(\frac{l-m}{N})\\
&+ \sum_{m=l'}^{l}\left( \sum_{u=0}^{l-m}
\beta^{\alpha}_{u} \overline {\beta^{\alpha}_{m+u}} \right)f(\frac{l-m}{N}).
\end{align*}
Then using the asymptotic expansion of the coefficients $\beta_{u}^{\alpha}$ 
and the hypotheses on $f$ we obtain 
\begin{enumerate}
\item
\begin{align*}
\Bigl \vert 
& \sum_{m=0}^{l'}\left( \sum_{u=0}^{l'-m}
\beta^{\alpha}_{u} \overline {\beta^{\alpha}_{m+u}} \right)
\left( f(\frac{l-m}{N}) -  f(\frac{l'-m}{N})\right)\Bigr \vert\\
& \le K_{1} \frac{l-l'}{N} N^{2\alpha}
\left( \frac{1}{\Gamma^{2}(\alpha)} \int _{0}^{l'/N}
 \int_{0}^{l'/N-z} t^{\alpha-1}(z+t)^{\alpha-1}dt dz +o(1)\right)\\
 & \le K_{2} \frac{l-l'}{N}N^{2\alpha}
\left(\frac{1} {\Gamma^{2}(\alpha)} \int _{0}^{1}
 \int_{0}^{1} t^{\alpha-1}(z+t)^{\alpha-1}dt dz +o(1)\right)
 \end{align*}
 \item
 \begin{align*}
\Bigl \vert 
 \sum_{m=0}^{l'}\left( \sum_{u=l'-m}^{l-m}
\beta^{\alpha}_{u} \overline {\beta^{\alpha}_{m+u}} \right)
 f(\frac{l-m}{N}) )\Bigr \vert
& \le K_{3} N^{2\alpha-1} \frac{(l-l')}{N}  \sum_{m=0}^{l^{\prime}} (\frac{l-m}{N})^{\alpha-1}
 (\frac{l'}{N})^{\alpha-1} f(\frac{l-m}{N}) \\
& \le K_{4} N^{2\alpha}\frac{l-l'}{N} \Vert f\Vert _{\infty}
 \end{align*}
 \item
 Lastly since the function $f$ is contracting 
 we obtain 
  \begin{align*}
\Bigl \vert 
 \sum_{m=l'}^{l}\left( \sum_{u=0}^{l-m}
\beta^{\alpha}_{u} \overline {\beta^{\alpha}_{m+u}} \right)
 f(\frac{l-m}{N}) )\Bigr \vert
&\le K_{5} \vert \frac{l-l'}{N} \vert  \sum_{m=l'}^{l} \sum_{u=0}^{l-m}\beta_{u}
 \overline {\beta^{\alpha}_{m+u}} \vert f \Vert_{\infty} \\
 &\le K_{6} \vert \frac{l-l'}{N} \vert N^{2\alpha}
 \end{align*}
 \end{enumerate}
 with $K_{1}, \cdots, K_{6}$ which are not dependent on $N$.
 So we can conclude that we have  a  positive constant  $M$
 such that   
 $ \vert \psi_{1,1,N} (l) -\psi_{1,1,N} (l')\vert \le M \frac{l-l'}{N} N^{2\alpha} $
 for $\frac{l}{N}, \frac{l'}{N} \in [\delta _{1}, \delta _{2}]$.
 To get the same result for $\psi_{1,2,N}$ 
we introduce decomposition 
\begin{align*}
\psi_{1,2,N}(l) &= \sum_{m=0}^{N-l'}  \sum_{v=0}^{l'} \beta_{v}^{\alpha} 
\overline{\beta_{m+v}^{\alpha}} f(\frac{m-l}{N})
- \sum_{m=N-l}^{N-l'}  \sum_{v=0}^{l'} \beta_{v}^{\alpha} 
\overline{\beta_{m+v}^{\alpha}} f(\frac{m-l}{N})\\
&+ \sum_{m=0}^{N-l}  \sum_{v=l'+1}^{l} \beta_{v}^{\alpha} 
\overline{\beta_{m+v}^{\alpha}} f(\frac{m-l}{N})
\end{align*}
Then the same methods as above allow us to conclude.

  Let's now show  the lemma for the difference 
  $\vert \psi_{2,N}(l) - \psi_{2,N}(l')\vert$.
  As previously we can split $\psi_{2,N}(l)$ in $\psi_{2,1,N}(l)+\psi_{2,2,N}(l)$
  with 

$$\psi_{2,1,N}(l) = \sum_{m=0}^l \left( \sum_{u=0}^m\frac{ \beta^\alpha_{m-u}}{N}
  \sum_{v=0}^{l-u} \overline{\beta^\alpha_{l-u-v}F_{N,\alpha} 
  (\frac{v}{N})}\right) f(\frac{m}{N})$$
  and
  $$\psi_{2,2,N}(l) = \sum_{m=l}^N \left( \sum_{u=0}^l\frac{ \beta^\alpha_{l-u}}{N}
  \sum_{v=0}^{m-u} \overline{\beta^\alpha_{l-u-v}F_{N,\alpha} 
  (\frac{v}{N})}\right) f(\frac{m}{N}).$$
  First $\psi_{2,1,N}(l)$ is also 
 $$\psi_{2,1,N}(l) = \sum_{m=0}^l \left(\sum_{u=0}^{l-m}\frac{ \beta^\alpha_{u}}{N}
  \sum_{v=0}^{m+u} \overline{\beta^\alpha_{v}F_{N,\alpha} 
  (\frac{m+u-v}{N})}\right) f(\frac{l-m}{N}).$$
Finally we have (always for $l'<l$)
\begin{align*}
\psi_{2,1,N}(l) &= \sum_{m=0}^{l'}  f(\frac{l-m}{N}) \left( \sum_{u=0}^{l'-m} 
  \frac{\beta^\alpha_{u}}{N} 
  \sum_{v=0}^{m+u} \overline{\beta_{v}^{(\alpha)} F_{N,\alpha}(\frac{m+u-v}{N})}\right)\\
  &+  \sum_{m=0}^{l'} f(\frac{l-m}{N}) \sum_{u=l'-m}^{l-m}
\frac{\beta^{\alpha}_{u}}{N}
 \left( \sum_{v=0}^{m+u}\overline{ \beta^\alpha_{v}
  F_{N,\alpha} (\frac{m+u-v}{N})} \right)\\
  &+ \sum_{m=l'}^l  f(\frac{l-m}{N}) \sum_{u=0}^{l-m}
\frac{\beta_{u}^{\alpha} }{N} \sum_{v=0}^{m+u} 
\overline{\beta_{v} ^{\alpha}  F_{N,\alpha} (\frac{m+u-v}{N})}
\end{align*}
Hence to bound the difference $\vert \psi_{2,1,N}(l) -\psi_{2,1,N}(l')\vert $
with $l'<l$ 
   we have three terms to consider.
   \begin{enumerate}
   \item
   First 
  $$ \sum_{m=0}^{l'} \left( f(\frac{l-m}{N}) - f(\frac{l'-m}{N})\right)\left( \sum_{u=0}^{l'-m} 
  \frac{\beta_{u}^\alpha}{N} 
  \sum_{v=0}^{m+u}\overline{ \beta_{v}^{\alpha} F_{N,\alpha}(\frac{m+u-v}{N})}\right)
$$
using (\ref{eqF})  we can bound this term by 
$$ K_{7}\frac{l-l'}{N} N^{2\alpha} \frac{1}{N} \sum_{m=0}^{l'} 
\frac{(l'-m)^{\alpha}}{N} \left( \frac{l'}{N}\right)^{\alpha} 
\vert \ln (1-\delta_{2})\vert = K_{8}\frac{l-l'}{N} N^{2\alpha}.$$
\item
 Then we have to study
$$ \sum_{m=0}^{l'} f(\frac{l-m}{N}) \sum_{u=l'-m}^{l-m}
\frac{\beta_{u}}{N}
 \left( \sum_{v=0}^{m+u}\overline{ \beta^\alpha_{v}
  F_{N,\alpha} (\frac{m+u-v}{N}) }\right) $$
 Always with (\ref{eqF}) we can bound this quantity by 
    $$ K_{9}\frac{l-l'}{N} \sum_{m=0}^{l'} f(\frac{l-m}{N}) N^{2\alpha-1}
  \left( \frac{l'-m}{N}\right)^{\alpha-1} 
  \left( \frac{l}{N}\right)^{\alpha} 
  \vert \ln (1-\delta_{2}) \vert = K_{10} \frac{l-l'}{N} N^{2\alpha}$$
   Lastly we have to consider
   $$\sum_{m=l'}^l  f(\frac{l-m}{N}) \sum_{u=0}^{l-m}
\frac{\beta_{u}^{\alpha} }{N} \sum_{v=0}^{m+u} 
\overline{\beta_{v} ^{\alpha}  F_{N,\alpha} (\frac{m+u-v}{N})}. $$
Always with the additional remark that $f$ is contracting on 
$[0,1]$ we can bound this quantity by 
$$K_{11} (\frac{(l-l')}{N})^2  N^{2\alpha} 
\left(\int_{0}{1} t^{\alpha-1} dt  \right)^2
 \vert \ln (1-\delta_{2}) \vert = K_{12} N^{2\alpha} (\frac{(l-l')}{N}).$$
 \end{enumerate}
 Since the quantities $K_{7}, \cdots, K_{12}$ are not dependent on $N$ we have got the property for this case.
   To treat the difference 
   $\vert \psi_{2,2,N}(l) - \psi_{2,2,N}(l')\vert $ we write, always for $l'<l$: 
   \begin{align*}
\psi_{2,2,N}(l) &= \sum_{m=0}^{N-l}  f(\frac{m-l}{N}) \left( \sum_{u=0}^{l'} 
  \frac{\beta^\alpha_{u}}{N} 
  \sum_{v=0}^{m+u}\overline{ \beta_{m+u-v}^{(\alpha)} F_{N,\alpha}(\frac{v}{N})}\right)\\
  &-  \sum_{m=N-l}^{N-l'}f(\frac{m-l'}{N}) \left( \sum_{u=0}^{l'} 
\sum_{v=0}^{m+u}\overline{\beta^\alpha_{m+u-v}
  F_{N,\alpha} (\frac{v}{N})} \right)\\
  &+ \sum_{m=0}^{N-l}  f(\frac{m-l'}{N}) \left(\sum_{u=l'}^{l}
\frac{\beta_{u}^{\alpha} }{N} \sum_{v=0}^{N-m+u} 
\overline{\beta_{v} ^{\alpha}  F_{N,\alpha} (\frac{v}{N})}\right)
\end{align*}
And the same methods as for $\psi_{2,1,N}$ provides that 
\begin{align*}
\Bigl \vert & \sum_{m=l}^N  \left( \sum_{u=0}^l \frac{\beta_{u}}{N} 
  \sum_{v=0}^{m-l+u} \overline\beta_{m-l+u-v}F_{N,\alpha} 
  (\frac{v}{N})\right) f(\frac{m}{N})\\
  -  &  \sum_{m=l'}^{N} \left(\sum_{u=0}^{l'} \frac{\beta_{u}}{N}
  \sum_{v=0}^{m-l'+u} \overline{\beta_{m-l+u-v}F_{N,\alpha} 
  (\frac{v}{N})}\right) f(\frac{m}{N})\Bigr \vert \le O( N^{2\alpha} \frac{l-l'}{N})
\end{align*}
For end the proof of the lemma we have to get the lemma for the difference 
$ \vert \psi_{3,N}(l)- \psi_{3,N}(l')\vert$ 
and $ \vert \psi_{4,N}(l)- \psi_{4,N}(l')\vert$.
These are the same ideas and methods used in the two proofs of the inequalities   
$ \vert \psi_{1,N}(l)- \psi_{1,N}(l')\vert$ 
and $ \vert \psi_{2,N}(l)- \psi_{2,N}(l')\vert$.
These same methods also make it possible to treat  the function 
$l \mapsto \displaystyle{\sum_{m=0}^{N}\sum_{u=0}^{\min(l,m)} \gamma_{N-\min (l,m)-u}^{\alpha} \overline \gamma_{N-\max(l,m)-u}}$.
\section{ Declarations}
 \subsection{Availability of data and material }
\textbf{Not applicable}
 \subsection{Competing interests}
\textbf{Not applicable}
 \subsection{Funding}
\textbf{Not applicable}
 \subsection{Authors' contributions}
\textbf{Not applicable}
\subsection{Acknowledgements}
\textbf{Not applicable}

  \bibliography{Toeplitzdeux}
\end{document}